\documentclass[12pt]{article}
\usepackage{amssymb}
\usepackage{amsmath}
\font\fFt=eusm10 
\font\fFa=eusm7  
\font\fFp=eusm5  
\def\K{\mathchoice
{\hbox{\,\fFt K}}
{\hbox{\,\fFt K}}
{\hbox{\,\fFa K}}
{\hbox{\,\fFp K}}}

\def\E{\mathchoice
{\hbox{\,\fFt E}}
{\hbox{\,\fFt E}}
{\hbox{\,\fFa E}}
{\hbox{\,\fFp e}}}

\def\etK{\eta_K}
\def\varK{{\varphi_K}}
\def\mur{{\mu(r)}}

\def\vark1{{\varphi_{1/K}}}
\def\rp{r^{{\prime}\thinspace 2}}

\def\var{\varphi}
\def\d{\delta}

\def\la{\lambda}

\def\R{\mathbb R}

\def\ds{\displaystyle}
\def\ce{{\E}}        
\def\kpr{{\K}'(r)}
\def\ck{{\K}}
\def\bs{\bigskip}

\title{Modular Equations and Distortion Functions}
\author{G. D. Anderson, S.-L. Qiu, and M. Vuorinen}
\date{}

\begin{document}

\numberwithin{equation}{section}

\maketitle

\begin{abstract}
\noindent
Modular equations occur in number theory, but it is less known
that such equations also occur in the study of deformation
properties of quasiconformal mappings.
The authors study two important plane
quasiconformal distortion functions,
obtaining monotonicity and convexity properties, and finding sharp
bounds for them.  Applications are provided that relate to the
quasiconformal Schwarz Lemma and to Schottky's Theorem. These
results also yield new bounds for singular values of complete
elliptic integrals.
\end{abstract}

{\bf Keywords and Phrases:}  Gr\"otzsch ring, elliptic integral, distortion
function, quasiconformal, quasiregular, quasisymmetric, quasi-M\"obius,
Schottky's Theorem, monotonicity, convexity.

{\bf 2000 Mathematics Subject Classification:}  Primary 30C62; Secondary
33E05.
\medskip

\setcounter{equation}{0}
\section{Introduction.}

For $r \in (0,1),$ let $\mu(r)$ denote the modulus of the plane
Gr\"otzsch ring $B\backslash [0,r]$,  where $B$ is the unit disk.  Then
\cite[p. 60]{20}
\begin{equation}
\mu(r) = \frac{\pi}{2} ~ \frac{\kpr}{{\K}(r)},
\end{equation}
where
\begin{equation}
\ck = \ck(r) \equiv \int^{\pi/2}_0\! (1\!-\!r^2 \sin^2 t)^{-1/2} dt,\;
\ck' = \kpr \equiv \ck(r'),
\end{equation}

\noindent $r' = \sqrt{1-r^2}, ~ 0 < r < 1$, are {\it complete elliptic
integrals
of the first kind} \cite{11}, \cite{12}, \cite{13}, \cite{40}.
In the sequel, we shall also need
the {\it complete elliptic integrals of the second kind}
\begin{equation}
\ce = \ce(r) \equiv \int^{\pi/2}_0 (1-r^2 \sin^2 t)^{1/2} dt, \quad
\ce' = \ce'(r) \equiv \ce(r'),
\end{equation}
for $r \in [0,1], ~ r'= \sqrt{1-r^2}.$

\normalsize

The function $\mu(r)$ and the distortion functions
$$
\begin{cases}
\varK (r) \equiv \mu^{-1}( \mur/K), \\[.2cm]
\varK (0) = 1 - \varK(1) = 0
\end{cases} \eqno{(1.4)}
$$
for $r \in (0,1)$ and $K \in (0, \infty)$, and
\setcounter{equation}{4}
\begin{equation}\label{ladef}
\la(K) \equiv \left[\frac{\varK(1/\sqrt{2})}{\var_{1/K} (1/\sqrt{2})}
\right]^2 = \frac{1- u^2}{u^2 } \, , u = \var_{1/K} (1/\sqrt{2}),
\end{equation}

\noindent $K \in (0, \infty)$, play an important role in the study of
deformation properties of quasiconformal maps
\cite{2}, \cite{7}, \cite{16}, \cite{18}, \cite{20}, \cite{22}, \cite{26},
\cite{38}. Note that here the identity (1.22) for $\varphi_K$ was used.
These functions have also found
applications in some other mathematical fields such as number theory
\cite{6}, \cite{9}, \cite{10}, \cite{11}, \cite{24}, \cite{37}, \cite{39}.

If  $R$ is a rectangle with sides parallel to the $x$- and $y$-axes of length $a$ and $b$, respectively, the {\it modulus} of $R$ is defined as $M(R) = b/a$.  If $Q$ is a (topological) quadrilateral, its {\it modulus} is defined as $M(Q) = M(R)$, where the rectangle $R$ is a conformal image of $Q$ and has sides parallel to the coordinate axes.  The modulus of a quadrilateral is a conformal invariant.

  Now let $f$ be a diffeomorphism of a domain $G$ in the $z$-plane onto a domain $G'$ in the $w$-plane.  The {\it maximal dilatation} at a point $z\in G$ is defined as
$$
K(z)=\inf_{G_z} \sup_{\overline Q\subset G_z}\frac{M(Q')}{M(Q)},
$$
where  $G_z$ denotes any subdomain containing $z$, $Q$ any quadrilateral with $\overline Q \subset G_z$, and $Q'$ the image of $Q$ under $f$.  If $K(z) \le K, 1 \le K < \infty$, for each $z\in G$, then  $f$ is called a (differentiable) $K$-{\it quasiconformal mapping}.  The definition of quasiconformality can be extended to the case of a homeomorphism $f$ that may fail to be differentiable at some points.  Next, for each $z\in G$ with $z\ne \infty \ne f(z)$, let
$$
H(z) = \limsup_{r\to 0}\frac{\max_{\varphi} |f(z+re^{i\varphi})-f(z)|}{\min_{\varphi}|f(z+ re^{i\varphi})-f(z)|}.
$$
If $z$ or $f(z)$ is infinite, the definition of $H(z)$ can be modified in an obvious way by means of inversion.  If  $f$  is continuously differentiable and has a non-zero Jacobian at $z$, then $H(z) = K(z)$ (see \cite{21}).  In \cite{21}, Lehto, Virtanen, and V\"ais\"al\"a proved that always  $H(z) \le \lambda (K(z))$, where $\lambda (K)$ is given by (\ref{ladef}), and this is the best possible result of this type.  
Moreover,  the boundary values of a $K$-quasiconformal automorphism $f$ of the upper
half plane satisfy the inequality
$$
\frac{1}{\la(K)} \le \frac{f(x+t)-f(x)}{f(x)-f(x-t)} \le {\la(K)}
 $$
for all real $x$ and $t, t\neq 0$, see \cite[p. 81]{20}. Accordingly,
several results in distortion theory of quasiconformal mappings
depend on $\la(K).$  One of the key areas of $K$-quasiconformal theory is the study of what happens when $K\to 1$.  Quantitative study of this subject requires explicit estimates for special functions, such as we are providing in this paper.

A {\it modular equation of degree} p is defined as
$$
\frac{\K'(s)}{\K(s)} = p \frac{\K'(r)}{\K(r)} ,
$$
with $s,r \in (0,1)$. Its solution is $s= \varphi_{1/p}(r) .$
Ramanujan introduced the convenient notation
$$
\alpha = r^2,\ \beta = s^2
$$
for use in connection with modular equations. With this notation,
a classical third-degree modular equation \cite[p. 105, (4.1.16)]{11}
reads as
$$
(\alpha \beta)^{1/4} +((1-\alpha)(1-\beta))^{1/4} = 1,
$$
with $\alpha = r^2$, $\beta = \varphi_{1/3}(r)^2$.
Ramanujan found numerous algebraic identities satisfied by
$\alpha$ and $\beta \,$ for several
prime numbers $p\,,$ see \cite{9}. Ramanujan also formulated, in unpublished
notes without proofs, many generalized modular equations.
Proofs were published in the landmark paper \cite{10} by Berndt, Bhargava,
and Garvin.

For each  positive integer $p$, there exists a unique number $k_p \in (0,1)$
such that
$$
\mu (k_p) = \frac{\pi \K '(k_p)}{2\K (k_p)} = \frac{\pi }{2} \sqrt{p}.
$$

\noindent The number $k_p$ is algebraic and is called the $p$th {\sl singular
value\/} of $\K (r)$ \cite [pp.\ 139, 296]{11}. Since $\mu (1/\sqrt{2})
= \pi /2$ we see that
$$
k_p = \mu ^{-1}(\sqrt{p} \mu (1/\sqrt{2})) = \varphi_{1/\sqrt{p}}
(1/\sqrt{2}).
$$

The many well-known facts about $\varphi_K$ yield information about $k_p$.
For instance, the infinite product expansion in \cite[Theorem 5.48 (3)]{7}
and the inequalities in \cite[Theorem 10.9 (1)]{7}, \cite{29}, or \cite{30}
imply such results.
Note that it follows from (1.5) that $k_p= 1/\sqrt{1+ \lambda(\sqrt{p})}\,.$
Many singular values have been found explicitly, because they have a
significant role in number theory. The
algebraic numbers $k_p , p = 1, \ldots, 9,$ are given in \cite [p.\ 139]{11}.
The values of $\K (k_p), p = 1, \ldots, 16,$ appear
in \cite [p.\ 298]{11}; from these one also obtains
$\K '(k_p), p = 1, \ldots, 16.$

In 1968, S. Agard \cite{1} introduced the following generalization of
$\la(K)$:
\begin{equation}
\etK(t) = \left[\frac{\varK(r)}{\var_{1/K} (r')} \right]^2, ~~ r =
\sqrt{\frac{t}{1+t}}, ~~ r' = \frac{1}{\sqrt{1+t}},
\end{equation}

\noindent for $K,t \in (0,\infty)$.  He showed that
\begin{equation}
\etK(t) = \sup_{f \in {\cal F}} \{|f(z)|:|z| = t\},
\end{equation}
where ${\cal F} = \{f:f$ is a $K$-quasiconformal
automorphism of the plane $\R^2$ with $f(0) = f(1) - 1 = 0 \}, K
\ge 1.$  Clearly $\etK(1) = \la(K)$.  
The function $\etK(t)$, called
the $\eta$-{\it distortion function}, has many important applications
in the study of quasiconformal, quasiregular, quasisymmetric, and
quasi-M\"obius mappings, and M\"obius-invariant metrics \cite{1}, \cite{6}, \cite{20}, \cite{32},
\cite{33}, \cite{34}, \cite{35}, \cite{36}.  Recently, G. Martin proved \cite[Theorem 1.1]{23}
that
\begin{equation}
\etK(t) = \sup \{|f(z)|:f \in {\cal A}(t), |z|=(K-1)/(K+1)\},
\end{equation}

\noindent where $t > 0$ and
$$
{\cal A}(t) = \{f: B \to\R^2 \backslash \{0,1\}\ |\ f ~~ \mbox{is
analytic with} ~~ |f(0)| = t\} \,.
$$
Hence, in the notation used in Schottky's Theorem \cite[p. 702]{15},
\begin{equation}
|f(z)| \le \etK(a) \equiv \Psi (a,|z|), ~ K = (1 + |z|)/(1 - |z|)\,,
\end{equation}

\noindent for $f \in {\cal A}(a), ~ a > 0,$ and $z \in B$.  Upper
bounds have been obtained for the function $\Psi(a,r)$ by W. K.
Hayman \cite{14}, J. Jenkins \cite{19}, J. Hempel \cite{17}, and S.
Zhang \cite{43}, while Zhang's estimates were recently improved in
\cite{27}.

In the past few years, many properties have been derived for the
special functions $\varK(r), ~ \la(K)$, and $\etK(t)$ \cite{3},
\cite{4}, \cite{5}, \cite{6}, \cite{7}, \cite{23}, \cite{25},
\cite{27}, \cite{29}, \cite{30}, \cite{31}, \cite{32},
\cite{37}, \cite{39}, \cite{41}, \cite{42}.  However, the
above-mentioned applications motivate us to study these functions
further.  For the applications, for example, we require better
estimates and some other properties for $\etK(t)$.

\bs

Since $\etK(t)$ is a generalization of $\la(K)$, it is natural
to ask how to extend the known results for $\la(K)$ to $\etK(t)$.
For example, $\la(K)$ has the following asymptotic expansion
\cite[Exercise 10.41 (6)]{7} (cf. \cite[(6.10), p. 82]{20}):
$$
\la(K) = \frac{1}{16} e^{\pi K} - \frac{1}{2} + \d(K), ~~ \mbox{where}
~~ 1 < e^{\pi K} \d(K) < 2, \eqno(1.10)
$$
for $K > 1$.  In Theorem 4.7 and Corollary 4.16 below we provide an
analogue of (1.10) for the function $\etK(t)$.
 \bs

In the present paper, the authors study the monotonicity and
convexity properties of certain combinations of the functions
$\la(K)$ and $\etK(t)$, from which sharp bounds for these functions
follow.  Note that, by (1.9), our sharp bounds for $\etK(t)$ also
give sharp estimates for the Schottky upper bound $\Psi(t,r)$.  Our
main results include the following Theorem 1.11, which leads to an
approximation of $\lambda(K)$ by a finite Taylor series, and Theorem
1.14, in which we extend to $\etK(t)$ the properties of $\la(K)$
proved in \cite[Lemma 3.54]{32}. Recalling that $k_p= 1/\sqrt{1+
\lambda(\sqrt{p})}\,, p=1,2,3,...,$ we see that the bounds for
$\lambda(K)$ in the next theorem also give bounds for the singular
values $k_p$. As far as we know, the resulting bounds for singular
values are new.

\bs

{\bf 1.11. Theorem.}  {\it Let $a = (4/\pi) \ck(1/\sqrt{2})^2 =
4.3768\dots$ and $c = a[4(a-1)^2 - a^2]/16 = 7.2372\dots$.
Then the function}
$$
f(K) \equiv (K-1)^{-3} \left[\la(K) - 1 - a(K-1) - \frac{1}{2}
a(a-1)(K-1)^2\right]
$$
{\it is strictly increasing from $(1, \infty)$ onto $(c,
\infty)$.  In particular, for} $K > 1$,
$$
\la(K) > 1 + a(K-1) + \frac{1}{2} a(a-1)(K-1)^2 + c(K-1)^3
\eqno(1.12)
$$
{\it and, for} $K \in (1,2)$,
$$
\begin{cases}
 1 + a(K-1) + \frac{1}{2} a(a-1)(K-1)^2 + c(K-1)^3 < \la(K) \\[.5cm]
\qquad \qquad < 1 + a(K-1) + \frac{1}{2} a(a-1) (K-1)^2 +
c_1(K-1)^3, \end{cases} \eqno(1.13)
$$
 {\it where}\ $c_1 \!=\! \la(2) - 1 - a - \frac{1}{2}
a(a-1) =\! 4 \sqrt{2} (\sqrt{2}+1)^2 - 1- a - \frac{1}{2} a(a-1)
=\! 20.2035\dots$.

\bs

{\bf 1.14. Theorem.}  {\it For $t \in (0, \infty)$, let $r =
\sqrt{t/(1+t)}$.}

(1)  {\it The function}
$$
f(K) \equiv \etK(t) e^{-2K \mu(r')}
$$
{\it is strictly increasing from $[1, \infty)$ onto}  $[te^{-2\mu(r')},
1/16).$

(2)  {\it The function}
$$
g(K) \equiv \frac{1}{\etK(t)} \frac{\partial \etK(t)}{\partial K}
$$
{\it is strictly decreasing from $(1, \infty)$ onto}  $(2\mu(r'), ~
(4/\pi) \ck(r) \ck'(r))$. {\it  In particular, for all $t \in (0,
\infty)$ and} $K \in (1, \infty),$
$$
te^{2(K-1) \mu(r')} < \etK(t) < te^{4(K-1) \ck(r) \ck'(r)/\pi}.
\eqno(1.15)
$$
(3)  {\it The function}
$$
h(K) \equiv  e^{-2 K \mu(r')}  \frac{\partial \etK(t)}{\partial K}
$$
{\it is strictly increasing from $(1, \infty)$ onto}  $((4/\pi) t
\ck(r) \ck'(r) e^{-2 \mu(r')}, \mu(r')/8)$.  {\it In particular, for
all $t \in (0, \infty)$ and} $K \in (1, \infty),$
$$
t + t\left(\frac{2}{\pi} \ck'(r)\right)^2 [e^{2(K-1)\mu(r')}-1] <
\etK(t)
< t + \frac{1}{16} e^{2 \mu(r')} [e^{2(K-1) \mu(r')}-1]\eqno{(1.16)}
$$
{\it and}
$$
1 + c(e^{\pi K} -e^{\pi}) < \la(K) < 1 + \frac{1}{16}(e^{\pi K}
-e^{\pi}), \eqno(1.17)
$$
{\it where} $c = \left( \frac{2}{\pi} \ck(1/\sqrt{2})\right)^2 e^{-\pi} =
0.0602 \dots$.
\bs

Throughout this paper, we let $r' = \sqrt{1-r^2}$ for $r \in [0,1]$.

We shall frequently employ the following well-known formulas
\cite[Appendix E]{7}, \cite{11}, \cite{16}, \cite{20}, \cite{40},
for $0 < r < 1, ~~ 0 <  K < \infty:$
$$
\frac{d \ck}{dr} = \frac{\ce - \rp \ck}{r \rp}, ~~ \frac{d
\ce}{dr} =
\frac{\ce - \ck}{r}, \eqno(1.18)
$$
$$
\ck \ce' + \ck' \ce - \ck\ck' = \frac{\pi}{2}, \eqno(1.19)
$$
$$
\frac{d \mur}{dr} = - \frac{\pi^2}{4 r \rp \ck(r)^2},
\eqno(1.20)
$$
$$
\begin{cases}
\ds{\dfrac{\partial s}{\partial r} = \dfrac{s}{Kr} \left(
\frac{s'\ck(s)}{r'\ck(r)}\right)^2 = \frac{ss^{{\prime} \thinspace 2}
\ck(s)
\ck'(s)}{r \rp \ck (r) \ck'(r)}}, \\[.4cm]

\ds{\frac{\partial s}{\partial K} = \frac{2}{\pi K}
ss^{{\prime}\thinspace 2} \ck (s)\ck'(s),}
\end{cases} \eqno(1.21)
$$
where $s= \varK(r),$
$$
\varK(r)^2 + \vark1(r')^2 = 1, ~~ \var_2(r) = \frac{2 \sqrt{r}}{1 +
r}. \eqno(1.22)
$$
We denote
$$
m(r) = \frac{2}{\pi}\rp \ck(r) \ck'(r), \eqno(1.23)
$$
for $r \in (0,1)$.  Then, by differentiation and (1.19), we have
$$
\frac{dm(r)}{dr} = \frac{1}{r} - \frac{4}{\pi r} \ce'(r) \ck(r).
\eqno(1.24)
$$

\bs

\setcounter{equation}{0}
\section{Preliminary results.}

In this section we prove two lemmas that will be needed for the proofs
of the main theorems in Sections 3 and 4.

\bs

{\bf 2.1. Lemma.}    (1)  {\it   The function $f(r) \equiv rr'\ck(r) \ck
(r')$ is strictly increasing on $(0,1/\sqrt{2}]$, and strictly decreasing on
$[1/\sqrt{2},1)$.}

(2)  {\it The function $g(r) \equiv \ck(r)/r$ is strictly decreasing on
$(0,1/\sqrt{2}]$.}

\bs

{\bf Proof.}  (1)  Differentiation gives
$$
r'f'(r) = \ck' \ce - \ck \ce',
$$
which is strictly decreasing from $(0,1)$ onto $(-\infty,\infty)$ and
has a unique zero at $r = 1/\sqrt{2}$.  Hence the result for $f$ follows.

(2)  By differentiation we obtain
$$
(rr')^2g'(r) =  (\ce - \rp \ck) - \rp \ck \equiv g_1(r),
$$
which is strictly increasing from $(0,1)$ onto $(-\pi/2,1)$ by
\cite[Theorem 3.21 (1), (7)]{7}, with $g_1(1/\sqrt{2}) = \ce(1/\sqrt{2}) -
\ck(1/\sqrt{2}) < 0.$  Hence the result for $g$ follows.  $\qquad
\square$

\bs

{\bf 2.2. Lemma.}    (1)  {\it  The function $f(r) \equiv (2/\pi) \ck(r)
\ck'(r) + \log r - \mu(r')$ is strictly decreasing from $(0,1)$ onto}
$(0,\log 4)$.

(2)  {\it The function $g(r) \equiv (2/\pi) \ck(r) \ck'(r) + \log
(r'/r)$ is strictly decreasing from $(0,1)$ onto $(\log 4, \infty)$.}

\bs

{\bf Proof.}  (1)  Since $f(r)$ can be written as $f(r) = (m(r) + \log
r) + (m(r') + \log r') - (\mu(r') + \log r')$, it follows from
\cite[Theorem 3.30 (1)]{7} and \cite[(2.11)]{20}  that
$f(0^+) = \log 4$ and $f(1^-) = 0$.

 From (1.24)  and (1.19) we have
the formula
$$
\frac{d}{dr} (m(r) + \log r) = - \frac{4}{\pi r} \ck'(\ck - \ce).
\eqno(2.3)
$$
Then, differentiating and using (1.19), (1.20), and (2.3), we get
$$
\frac{\pi}{4} r(r'\ck')^2 f'(r) = f_1(r) \equiv
\ck^{{\prime}\thinspace 3} (\ce -
\rp \ck) - \frac{\pi}{4} \left(r^2 \ck^{{\prime}\thinspace 2}
+ \frac{\pi^2}{4}\right).
\eqno(2.4)
$$
It follows from \cite[Theorem 3.21 (1)]{7}  that
$$
\begin{cases}
\ds{f_1(r)  = r^2 \ck^{{\prime} \thinspace 2}\left(\ck' \frac{\ce-\rp
\ck}{r^2} - \frac{\pi}{4}\right) - \frac{\pi^3}{16}}  \\[.5cm]
\ds{\qquad \quad < (r \ck')^2 \left(\ck' - \frac{\pi}{4}\right) -
\frac{\pi^3}{16}
\equiv f_2(r).} \end{cases} \eqno(2.5)
$$
Differentiation gives
$$
\frac{\rp f'_2(r)}{r \ck^{{\prime} \thinspace 2}(\ck' - \ce')} =
f_3(r) \equiv  2 - \frac{\pi}{2 \ck'} - \frac{\ce' - r^2\ck'}{\ck' -
\ce'}.
$$
Since $(\ce' - r^2 \ck')/(\ck' - \ce')$ is strictly increasing from
$(0,1)$  onto $(0,1)$ \cite[Theorem 1.8 (2)]{28}, $f_3$ is strictly
decreasing from $(0,1)$ onto $(0,2)$, so that $f_2$ is strictly
increasing on $(0,1)$.  Clearly $f_2(1) = 0$.  Hence the monotonicity of
$f$ follows from (2.4) and (2.5).

\medskip

(2) By differentiation and (1.19), we get
$$
g'(r) = -\frac{4 \ck}{\pi r \rp} (\ce' - r^2 \ck'),
$$
which is negative by \cite[Theorem 3.21 (1)]{7}.  Hence the monotonicity of
$g$ follows.

Since $g(r) = (m(r) + \log r) + (m(r') + \log r') - 2 \log r,$  the
limiting values of $g$ follow from \cite[Theorem 3.30 (1)]{7} . $\qquad
\square$

\bs

\setcounter{equation}{0}
\section{Properties of $\pmb{\la(K)}$.}

In this section,  we prove several refinements of some known results
for the function $\la(K)$.  Our first result improves
\cite[Corollary 10.33, Theorem 10.35]{7}  and \cite[Corollary 2.11]{37}.

\bs

{\bf 3.1. Theorem.}  (1)  {\it The function $f(K) \equiv (\log
\la(K))/(K-1)$ is strictly decreasing and convex from $(1, \infty)$
onto $(\pi, a)$, where} $a = (4/\pi) \ck (1/ \sqrt{2})^2 = 4.3768 \dots$.

(2) {\it The function $g(K) \equiv (1/(K-1)) \log (\la(K)
e^{b(-K+1/K)})$, where $b = a/2$, is strictly increasing from $(1,
\infty)$ onto} $(0, \pi-b)$.

{\bf Proof.} Let $r = \mu^{-1} (\pi/(2K))$.  Then, by (1.1), (1.5), and (1.22),
$\la(K) = (r/r')^2$.  Using (1.20), we get
$$
\frac{dr}{dK} = \frac{2}{\pi} r \rp \ck' (r)^2, ~~
\frac{d \la(K)}{dK}= \frac{4}{\pi} \la (K) \ck'(r)^2. \eqno(3.2)
$$

(1)  By differentiation and (3.2), we have
$$
f'(K) = f_1(K)/f_2(K), \eqno(3.3)
$$
where $f_1(K) = \dfrac{4}{\pi}(K-1) \ck'(r)^2 - \log \la(K)$ and $f_2(K)
= (K-1)^2,$
$$
f'_1(K) = - \left( \frac{4}{\pi} \right)^2 (K-1) \ck'(r)^3 [\ce'(r)
- r^2 \ck'(r)], \eqno(3.4)
$$
and
$$
f'_1(K)/f'_2(K) = -  \frac{8}{\pi^2} \ck'(r)^3
[\ce'(r) - r^2 \ck'(r)] \equiv f_3(r).
$$

 From (3.4) and \cite[Theorem 3.21 (1)]{7}, we see that $f_1$ is strictly
decreasing on $(1, \infty)$ with $f_1(1) = 0$, so that, by (3.3), $f$
is strictly decreasing on $(1, \infty)$.

By \cite[Theorem 3.21 (1)]{7}, $f_3$ is strictly increasing on $(0,1)$, so
that $f'_1(K)/f'_2(K)$ is strictly increasing in $K$ on $(1, \infty)$,
and so is $f'$ by the l'H\^opital Monotone Rule \cite[Theorem 1.25]{7}.
This yields the convexity of $f$.

The limiting values of $f$ follow from \cite[Corollary 2.11]{37}.

\medskip

(2)  Let $g_1(K) = \log \la(K) + b(\frac{1}{K} - K)$ and $g_2(K) =
K-1$.  Then, by (3.2) and differentiation, since $K = \ck (r)/\ck'(r)$,
we have
$$
g'_1(K) /g'_2(K) = g_3(K) \equiv \frac{4}{\pi} \ck'(r)^2 - b(1 +
K^{-2}), \eqno(3.5)
$$
$$
\frac{1}{2} K^3 g'_3(K) = b - \frac{8}{\pi^2} \ck(r)^3[\ce'(r) - r^2
\ck' (r)] \equiv g_4(r). \eqno(3.6)
$$

It follows from \cite[Theorem 3.21 (1),(7)]{7}  that $g_4$ is strictly
increasing on $(0,1)$.  Since $r \in (1/\sqrt{2},1) $ and since, by
(1.19),
$$
\begin{cases}
g_4(1/\sqrt{2}) &= \frac{2}{\pi} b \left[ \frac{\pi}{2} - 2
\ck(1/\sqrt{2})\left(\ce (1/\sqrt{2}) - \frac{1}{2}
\ck(1/\sqrt{2})\right) \right]\\
&=\frac{2}{\pi} b\left[ \frac{\pi}{2} - 2 \ck
(1/\sqrt{2}) \ce(1/\sqrt{2}) + \ck(1/\sqrt{2})^2 \right] = 0,
\end{cases}
$$
it follows from (3.6) that $g_3$ is strictly increasing on $(1,\infty)$.
Hence the monoto\-nicity of $g$ follows from (3.5) and \cite[Theorem 1.25]{7}.

The limiting values of $g$ follow from l'H\^opital's Rule and (3.5).
$\qquad \square$

\bs

{\bf 3.7. Corollary.}  {\it For} $K > 1$,
$$
\max \{ e^{\pi(K-1)}, e^{b(K-1/K)} \} < \la(K) < \min
\{e^{a(K-1)}, e^{(\pi+b/K)(K-1)} \},\eqno(3.8)
$$
{\it where $a$ and $b$
are as in Theorem 3.1.  Moreover,}
$$
\lim_{K \to 1} \la(K)^{1/(K-1)} = e^a, ~~ \lim_{K \to \infty}
\la(K)^{1/K} = e^{\pi}. \eqno(3.9)
$$

\bs

{\bf Proof.}  The estimates in (3.8) follow immediately from Theorem
3.1.  The lower estimates and first upper estimate in (3.8) also
follow from \cite[Corollary 10.33, Theorem 10.35]{7},
while (3.9) follows directly from (3.8).  $\qquad \square$

\bs

In \cite[Lemma 12, p. 80]{8}, P. P. Belinski gave the inequality
$$
\la(K) < 1 + 12 (K-1) \eqno(3.10)
$$
for $K > 1$ close to 1.  However, his proof given in \cite[pp. 80-82]{8}
for (3.10) is not valid (cf. \cite[p. 412]{37}).  Corollary 3.5 of \cite{37}
gives an improved form of (3.10).  Theorem 1.11 is related to this
kind of property of $\la(K)$ for $K > 1$ close to 1, and improves
\cite[Corollary 3.5]{37}.

\bs

{\bf 3.11. Proof of Theorem 1.11.}  Let
$r = \mu^{-1} (\pi/(2K)), ~ \ck = \ck(r), ~ \ck' = \ck'(r)$,
and $\ce' = \ce' (r)$.  Then (3.2) holds, and by (1.18) and
(3.2) we have
$$
\begin{cases}
\hfill\dfrac{d\ck'}{dK} &=
-\dfrac{2}{\pi} \ck^{{\prime}\thinspace 2}(\ce'- r^2 \ck'),   \\
\vspace{.1in}\la''(K) &=
\dfrac{4}{\pi} \la'(K) \ck'[\ck' - (\ce' - r^2\ck')],\\
\vspace{.1in}\la'''(K) &=
\dfrac{6}{\pi} \la''(K) \ck'[\ck' - (\ce' - r^2 \ck')]
-\dfrac{8}{\pi^2} (1-r^2 \rp) \la'(K) \ck^{{\prime}\thinspace 4}.
\end{cases} \eqno(3.12)
$$

By (1.19), we have
$$
\ck(1/\sqrt{2})\left[\ce(1/\sqrt{2}) - \frac{1}{2}
\ck(1/\sqrt{2})\right] = \frac{\pi}{4}. \eqno(3.13)
$$
Using (3.2), (3.12), and (3.13), we obtain
$$
\la'(1) = a, ~~ \la''(1) = a(a-1), ~~ \la'''(1) = 6c. \eqno(3.14)
$$

Next, let $g(K) = \la(K) - 1 - a(K-1) - \frac{1}{2} a(a-1)(K-1)^2$ and
$h(K) = (K-1)^3$.  Then
$$
\begin{cases}
\ds{\frac{g'(K)}{h'(K)} = \frac{\la'(K) - a - a(a-1)(K-1)}{3(K-1)^2}},
\\[.6cm]
\ds{\frac{g''(K)}{h''(K)} = \frac{\la''(K) - a(a-1)}{6(K-1)}, ~~
\frac{g'''(K)}{h'''(K)} = \frac{1}{6} \la'''(K)}.
\end{cases} \eqno(3.15)
$$

 From (3.2), (3.12), and the fact that $\la(K) = (r/r')^2$, where $r
= \mu^{-1}(\pi/(2K))$, it follows that
$$
\la'''(K) = \frac{32}{\pi^3} (r \ck^{{\prime}\thinspace 2})^2 F(K),
\eqno(3.16)
$$
where
\begin{align*}
F(K) & = \frac{3}{\rp} [\ck' - (\ce'-r^2 \ck')]^2 - \frac{1-r^2
\rp}{\rp} \ck^{{\prime}\thinspace 2} \\
& = 2 \left(\frac{\ck'}{r'}\right)^2 + (r \ck')^2 - 3 \frac{\ce'-r^2
\ck'}{\rp} [(\ck'+ \ce') + (1+r^2) \ck'-2 \ce'].
\end{align*}
Using (3.13), we get
\begin{align*}
F(1) & = 6 \left\{\ck (1/\sqrt{2}) - \left[\ce (1/\sqrt{2}) - \frac{1}{2}
\ck(1/\sqrt{2})\right]\right\}^2 - \frac{3}{2} \ck(1/\sqrt{2})^2 \\
& = \frac{9}{2} \ck (1/\sqrt{2})^2 + 6\left[ \ce (1/\sqrt{2}) -
\frac{1}{2} \ck(1/\sqrt{2})\right]^2 - 3 \pi = 7.1210 \dots > 0.
\end{align*}
Since $r \in (1/\sqrt{2},1),$  it follows from Lemma 2.1 (2) and
\cite[Theorem 3.21(7)]{7},
 that $\ck'/r'$ and $r \ck'$ are increasing, while by
\cite[Theorem 3.21 (1), Lemma 1.33 (3)]{7}
$(\ce' - r^2 \ck')/\rp$ and $\ck' + \ce'$
are decreasing. Next, by \cite[Exercise 3.43 (4)(a)]{7},
$(1+r^2) \ck' - 2 \ce'$
is strictly decreasing from $(0,1)$ onto $(0, \infty)$.  Thus we
conclude that
$F(K)$ is strictly increasing in $K$ on $(1,\infty)$.  Hence, it follows
from \cite[Theorem 3.21 (7)]{7}  and (3.16) that $\la'''$ is
strictly increasing on  $(1, \infty)$.

By (3.14), we observe that $g'(1) = h'(1) = g''(1) = h''(1) = 0$.
Hence the monotonicity of $f$ follows from (3.15) and \cite[Theorem 1.25]{7}.

By l'H\^opital's Rule, (3.14), and (3.15), we get $f(1^+) =
c$.  Since $K = \ck/\ck', ~~ \la(K) = (r/r')^2,$ and
$$
f(K) = G(r) = \frac{r^2 - \rp \{1+a[(\ck/\ck')-1] + a_1[(\ck/\ck')
-1]^2\}}{\rp [(\ck/\ck') - 1]^3},
$$
where $a_1 = \frac{1}{2} a(a-1),$ and since $\sqrt{r'} \ck \to 0$ as $r
\to 1$ \cite[Theorem 3.21(7)]{7}, we see that
$$
\lim_{K \to \infty} f(K) = \lim_{r \to 1} G(r) = \infty.
$$

Finally, (1.12) is clear.  Since
$$
f(2) = \la(2) - 1 - a - \frac{1}{2} a(a-1) = c_1,
$$
(1.13) follows from the monotonicity of $f$. We have used
\cite[Theorem 10.5(4)]{7} in the evaluation of $\lambda (2)$. $\qquad \square$

\bs

{\bf 3.17. Corollary.}  (1) {\it As} $K \to 1$,
$$
\la(K) = 1 + a(K-1) + \frac{1}{2} a (a-1)(K-1)^2 + O((K-1)^3),
$$
{\it where $a$ is as in Theorem} 1.11.
\medskip

(2) {\it Let $\d > 0$ be an arbitrary real number and let $c_1$ be as
in Theorem} 1.11. {\it Then, for}
$1 < K \le K_0\equiv 1 + (\sqrt{a^2_1 + 4c_1 \d} ~- a_1)  /(2c_1),$
$$
\la(K) < 1 + (a + \d)(K-1). \eqno(3.18)
$$
{\it In particular, for}
$$
1 < K \le 1 + \left(\sqrt{a^2_1 + 4c_1(5-a)} ~
-a_1\right)/(2c_1) = 1.07066 \dots,
$$
$$
\la(K) < 1 + 5(K-1). \eqno(3.19)
$$

\bs

{\bf Proof.}  Part (1) follows immediately from (1.13).  By
(1.13), we see that (3.18) holds if
$$
c_1(K-1)^2 + a_1(K-1) - \d \le 0. \eqno(3.20)
$$
Clearly, for $K > 1$, (3.20) holds if and only if $K \le K_0$.

Taking $\d = 5-a$, we obtain (3.19) from (3.18). $\qquad \square$

\bs

The next result improves (1.10).

\bs

{\bf 3.21. Theorem.}  {\it For} $K > 1$,
$$
\frac{1}{16} e^{\pi K} - \frac{1}{2} + c_1(K) e^{-\pi K} < \la(K) <
\frac{1}{16} e^{\pi K} - \frac{1}{2} + c_2(K) e^{-\pi K}, \eqno(3.22)
$$
{\it where}
$$
c_1(K) \equiv \frac{1}{4} + \frac{1}{4e^{-2 \pi K}+1} > \frac{21}{16} -
c,
$$
$$
c_2(K) \equiv \frac{1}{16} \left(1 + 4 \frac{5e^{-4 \pi K} + 14e^{-2 \pi
K}+5}{e^{-6 \pi K} + 7e^{-4 \pi K} + 7e^{-2 \pi K} + 1} \right) <
\frac{21}{16},
$$
and $c \equiv (68 + e^{2 \pi})/[16(4 + e^{2 \pi})] = 0.06991 \dots$.

\bs

{\bf Proof.}  It follows from \cite[(5.20), (5.4),
Theorem 5.13 (4)]{7}  that, for all $r \in (0,1)$,
$$
\log \frac{1 + \sqrt[4]{r'}}{1 - \sqrt[4]{r'}} < 2 \mur < \log
\frac{2(1 + \sqrt{r'})}{1 - \sqrt{r'}}, \eqno(3.23)
$$
from which we get, with $x \equiv \exp(2 \mur) \ge 2$,
$$
1 - \left( \frac{x-1}{x+1}\right)^8 < r^2 < 1 - \left(
\frac{x-2}{x+2} \right)^4 \eqno(3.24)
$$
for all $r \in (0, \mu^{-1} (\log \sqrt{2})]$.

Let $r = \mu^{-1}(\pi K/2).$  Then $0 < r < 1/\sqrt{2} < \mu^{-1}(\log
\sqrt{2}), ~~ x = e^{\pi K},$ and $\la (K) = (r'/r)^2.$

It follows from the second inequality in (3.24) that
\begin{align*}
\la(K) & = \frac{1}{r^2} - 1 > \frac{(x+2)^4}{16 x(x^2+4)} - 1 =
\frac{e^{\pi K}}{16} ~ \frac{(1+2e^{-\pi K})^4}{1 + 4e^{-2 \pi K}} -1
\\[.3cm]
& = \frac{e^{\pi K}}{16} \left[1 + 8e^{-\pi K} + 4 e^{-2 \pi K}\left(1 +
\frac{4}{4e^{-2 \pi K} + 1} \right) \right] -1 \\[.3cm]
& = \frac{1}{16} e^{\pi K} - \frac{1}{2} + \frac{1}{4} e^{-\pi K}\left(
1 + \frac{4}{4e^{-2 \pi K} + 1} \right),
\end{align*}
and hence the first inequality in (3.22) holds.  Clearly,
$$
c_1(K) > \frac{1}{4} + \frac{1}{4e^{-2 \pi} + 1} = \frac{21}{16} -
c.
$$

Next, it follows from the first inequality in (3.24) that
\begin{align*}
\la(K) & = \frac{1}{r^2} - 1 < \frac{x (1+x^{-1})^8}{16
(1+x^{-2})(x^{-4} + 6x^{-2}+1)}- 1  \\[.3cm]
& = \frac{1}{16} e^{\pi K} + \frac{e^{\pi K}}{16} ~ \frac{(1+y)^8 -
(1+y^2)(1+6y^2 + y^4)}{(1+y^2)(1+6y^2 + y^4)} -1
\\[.3cm]
& = \frac{1}{16} e^{\pi K} + \frac{e^{\pi K}}{16} \left[8y + y^2 \left(
1 +  4  \frac{5 y^4 + 14 y^2 + 5}{y^6 + 7y^4 + 7y^2 + 1} \right)
\right] -1,
\end{align*}
where $y = 1/x = e^{- \pi K}$, and hence the second inequality in
(3.22) follows.

Finally, it is easy to verify that $(5t^2 + 14t + 5)/(t^3 + 7t^2 + 7t +
1)$ is a strictly decreasing function of $t$ on $(0, \infty)$.  Hence
$c_2(K)$ is strictly increasing on $[1, \infty)$, so that  $c_2(K) <
\lim_{K \to \infty} c_2(K) = 21/16. \qquad \square$

\bs

{\bf 3.25. Remark.}  Computation gives:
$$
\frac{1}{16} e^{\pi } - \frac{1}{2} + c_1(1) e^{-\pi}
 = 0.9999902 \dots,
$$
$$
\frac{1}{16} e^{\pi} - \frac{1}{2} + c_2(1) e^{-\pi} = 1.0025922 \dots .
$$
Hence, even when $K$ is close to 1, the lower and upper estimates for
$\la(K)$ given in (3.22) are very close to each other.

\bs

\setcounter{equation}{0}
\section{Properties of $\pmb{\etK(t)}$.}

In this section, we study some properties of $\etK(t)$.  We first
extend to $\etK(t)$ the properties of $\la(K)$ proved in \cite[Lemma 3.54]{32}.

\bs

{\bf 4.1. Proof of Theorem 1.14.}  Let  $s = \varK(r)$.
Then $\etK(t) = (s/s')^2$ and, using (1.4) and (1.21), we have
$$
\frac{\partial \etK(t)}{\partial K} = \frac{4}{\pi K} \etK(t) \ck
(s) \ck'(s) = \frac{2}{\mur} \ck'(s)^2 \etK(t). \eqno(4.2)
$$

(1) We may rewrite $f(K)$ as
$$
f(K) = \left[\frac{s}{s'e^{\mu(s')}} \right]^2 = e^{-2[\mu(s') +
\log(s'/s)]}.
$$
Hence the monotonicity of $f$ follows from \cite[Theorem 5.13 (3)]{7}.  The
limiting values are clear.

(2) By (4.2) we may write $g(K)$ as
$$
g(K)=\frac{2\ck'(s)^2}{\mu(r)},
$$
from which the monotonicity follows immediately.
We obtain (1.15) by integrating the inequalities
$$
2\mu(r')<g(K)<\frac {4}{\pi}\ck(r)\ck'(r)
$$
with respect to $K$ over $(1,K)$ and using the fact that $\eta_1(t)=t$.

(3)  By (4.2), we have
$$
h(K) = \frac{2}{\mur} \left[ \frac{\sqrt{s} ~
\ck'(s)}{(s'/\sqrt{s})e^{\mu(s')}} \right]^2.
$$
Hence the monotonicity of $h$ follows from \cite[Theorem 3.21(7)]{7}
and \cite[Lemma 3.54(1)]{32}.  We obtain (1.16) and (1.17) by integrating
the inequalities
$$
\frac{4}{\pi} t \ck (r)\ck'(r)e^{2(K-1)\mu(r')} <
\frac{\partial \etK}{\partial K} < \frac{1}{8} \mu(r') e^{2K\mu(r')}
$$
with respect to $K$ over $(1,K)$.  $\qquad \square$

\bs

{\bf 4.3. Remark.} An injective mapping $f:X\to Y$, where $X$ and
$Y$ are metric spaces with distances denoted by $|a-b|$, is called
{\it quasisymmetric} if there exists a homeomorphism
$\eta: [0,\infty )\to [0,\infty )$, $\eta(0) = 0$, such that for all
$a,b,c\in X$, $a\ne c$,
$$
\frac{|f(a)-f(b)|}{|f(a)-f(c)|} \le \eta \left(\frac{|a-b|}{|a-c|}\right).
$$
Further, if $s >0,$ we say that $f$ is $s$-{\it quasisymmetric} if $f$
is quasisymmetric with $\eta(t) \le t+s$ for all
$t \in (0,\max \{1, 1/s \}) \,.$
Using (1.16), one can show that, for $K$ sufficiently close to 1, a
$K$-quasiconformal map $f$ of $\R^{2}$ into $\R^{2}$ is
$s$-quasisymmetric with $s = (K-1)^{4/9}$. This result improves the
conclusion stated in \cite[Remark 2.17]{36}, where it
was indicated that $f$ is $s$-quasisymmetric with $s = \sqrt[4]{K-1}$.
In fact, by \cite[Theorem 2.16]{36} and (1.16), we only need to show
that, for $K$ sufficiently close to 1,
$$
e^{2\mu(r)}[e^{2(K-1)\mu(r)}-1]/s \leq 16, \eqno(4.4)
$$
where $r = \sqrt{s/(1+s)}$ and $s = (K-1)^{4/9}$. Since $K=1+s^{9/4}$
and $s = (r/r')^{2}$, it follows from \cite[(2.11), p. 62]{20},
\cite[Theorem 3.21 (7)]{7}, l'H\^{o}pital's Rule, and (1.20) that
\begin{eqnarray*}
\lefteqn{\lim_{K\to 1} e^{2\mu(r)}\left[e^{2(K-1)\mu(r)}-1\right]/s} \\[.2cm]
&& = \lim_{r\to 0} e^{2(\mu(r)+\log r)} \frac{e^{2(r/r')^{9/2}\mu(r)}-1}{r^{4}} = 16
\lim_{r\to 0} \frac{e^{2(r/r')^{9/2}\mu(r)}-1}{r^{4}} \\[.2cm]
&& = \frac{8}{\pi} \lim_{r\to 0} \sqrt{r} [9{\cal K}(r){\cal K}'(r)-\pi] = 0.
\end{eqnarray*}
Thus there exists $K_{0} \in (1,2)$ such that (4.4) holds for
$K \in (1,K_{0}]$.

\bs
The next result extends \cite[Theorem 2.3]{6} to the function
$\etK (t)$ for $K > 0.$

\bs
{\bf 4.5. Theorem.}  {\it For each fixed $t \in (0, \infty), ~ F(K)
\equiv \etK(t) + 1$  is  log-convex as a function of $K$ on $(0,
\infty)$, while $\etK(t)$ is log-concave there.  In particular, for
$t,K,L \in (0, \infty)$ and $p,q \in (0,1)$ with} $p + q = 1,$
$$
\etK(t)^p \eta_L(t)^q < \eta_{pK+qL} (t) < [\etK(t) +
1]^p [\eta_L(t) + 1]^q - 1. \eqno(4.6)
$$

\bs

{\bf Proof.}  First, we observe that (4.2)  also holds for $K\in(0,1),$
so that the monotonicity properties of $f,g,$ and $h$ in Theorem 1.14 are
valid for $K \in (0,\infty)$.  Hence the log-concavity of $\etK(t)$
follows from Theorem 1.14 (2).  Alternatively, the logarithmic
derivative of $\etK(t)$ is $(2/\mur) \ck'(s)^2$, which is strictly
decreasing in $K$.

Since $\etK(t) = (s/s')^2$, where $s = \varK(r)$ and $r =
\sqrt{t/(1+t)}$, we have
$$
F(K) = \etK(t) + 1 = (s')^{-2},
$$
and hence, by (1.21),
$$
\frac{F'(K)}{F(K)} = \frac{2}{\mur} s^2 \ck^{\prime} (s)^2,
$$
which is strictly increasing in $K$ on $(0,\infty)$ by
\cite[Theorem 3.21 (7)]{7}.  The log-convexity of $F$ now follows.

The remaining conclusions are clear.  $\qquad \square$

\bs

The next result is a generalization of Theorem 3.21 to the function
$\etK(t)$, and provides an analogue of (1.10) for $\etK(t)$.

\bs

{\bf 4.7. Theorem.}   (1) {\it For $K > 1$ and} $0 < t < [\mu^{-1}
(\frac{1}{2K} \log 2)]^{-2}-1$,
$$
\etK(t) < [\mu^{-1} (\tfrac{1}{2} \log 2)]^{-2} -1 < 0.000011.
\eqno(4.8)
$$

(2) {\it For $K>1$ and} $t \ge [\mu^{-1} (\frac{1}{2K} \log 2)]^{-2} -
1,$
$$
\etK(t) > \frac{1}{16} e^{2 K \mu(r')} - \frac{1}{2} + c_1(t,K)
e^{-2K \mu(r')}, \eqno(4.9)
$$
{\it where $r = \sqrt{t/(1+t)}$ and}
$$
c_1(t,K) = \frac{1}{4} + \frac{1}{4e^{-4K\mu(r')} + 1} > \frac{1}{4}
+ \frac{1}{4e^{-4 \mu(r')} +1} > \frac{9}{20}. \eqno(4.10)
$$

(3) {\it For $t \in (0, \infty)$ and} $K > 1$,
$$
\etK(t) < \frac{1}{16} e^{2K \mu(r')} - \frac{1}{2} + c_2
(t,K)e^{-2K \mu(r')}, \eqno(4.11)
$$
{\it where $r$ is as in part} (2) {\it and}
$$
\begin{cases}
\ds{c_2(t,K) ~= \frac{1}{16} \left[ 1 + 4 \frac{5e^{-8K \mu(r')} +
14e^{-4K \mu(r')} + 5}{e^{-12K \mu(r')} + 7e^{-8K\mu(r')} + 7e^{-4 K
\mu(r')} + 1} \right]} \\[.4cm]
\ds{\qquad \qquad  < \frac{5}{16} + \frac{1}{4e^{-2K \mu(r')} + 1} <
\frac{21}{16}}.
\end{cases} \eqno(4.12)
$$

\bs

{\bf Proof.}  Let $u = \vark1(r'), ~ r = \sqrt{t/(1+t)}$, and $x = e^{2
\mu(u)} = e^{2K\mu(r')}$.  Then, by (1.6) and (1.22),
$$
\etK(t) = u^{-2} - 1 \eqno(4.13)
$$
and, by (3.23),
$$
\frac{1 + \sqrt[4]{u'}}{1 - \sqrt[4]{u'}} < x < 2 \frac{1 +
\sqrt{u'}}{1-\sqrt{u'}}~        . \eqno(4.14)
$$

(1)  For $K > 1$ and $0 < t < [\mu^{-1}(\frac{1}{2K} \log 2)]^{-2} -
1$, we have
$$
u = \vark1 (r') > \mu^{-1}(K \mu (\mu^{-1} (\tfrac{1}{2K} \log 2 ))) =
\mu^{-1} (\frac{1}{2} \log 2)
$$
so that, by (4.13),
$$
\etK(t) < [\mu^{-1} (\tfrac{1}{2} \log 2)]^{-2} - 1 < 0.99999476^{-2}
- 1 < 0.000011.
$$

(2) Since $t \ge [\mu^{-1} (\frac{1}{2K} \log 2)]^{-2} -1$, we have
$$
x = e^{2K \mu(1/\sqrt{1+t})} \ge e^{2K\mu(\mu^{-1}(\frac{1}{2K} \log
2))} = 2.
$$
Hence, it follows from the second inequality in (4.14) that
$$
u^2 < 1 - \left( \frac{x-2}{x+2} \right)^ 4,
$$
so that, by (4.13),
$$
\etK(t) > \frac{(x+2)^4}{16 x(x^2 + 4)} - 1.
$$
Then, by the method used in the proof of Theorem 3.21 one can easily
obtain (4.9).  The inequalities in (4.10) are clear.

(3)  Clearly, $x > 1$ for all $t > 0$ and $K > 1$.  Hence, it follows
from the first inequality in (4.14) that
$$
u^2 > 1 - \left( \frac{x-1}{x+1} \right)^8.
$$
Now, using (4.13), one can prove (4.11) by the method used in the proof
of the second inequality in (3.22), without any difficulty.

Finally, let
$$
f(y) = 4 \frac{5 y^2 + 56y + 80}{y^3 + 28y^2 + 112 y + 64} -
\frac{4}{y+1} - \frac{3}{4}.
$$
Then $c_2(t,K)$ can be written as
$$
c_2(t,K) = \frac{1}{4} \left[1+ \frac{4}{b^2 + 1} +f(b^2) \right],
\eqno(4.15)
$$
where $b = 2/x = 2e^{-2K\mu(r')}$.  Clearly, $b \in (0,2)$.  It is easy
to show that the function $f$, which can be rewritten as
$$
f(y) = 4 ~ \frac{4 y^3 +33y^2 + 24y +16}{(y+1)(y^3 + 28y^2 + 112 y +
64)}  - \frac{3}{4},
$$
is strictly decreasing on $(0,4)$ and that $f(0) = 1/4$.  Hence the
inequalities in (4.12) follow from (4.15).  $\qquad \square$

\bs

As a consequence of Theorem 4.7, the following corollary gives an
asymptotic expansion for $\etK(t)$ as $K \to \infty$ or $t \to \infty$.

\bs

{\bf 4.16. Corollary.}  {\it As $K \to \infty$ or $t \to \infty$,}
$$
\etK(t) = \frac{1}{16} e^{2K\mu(r')} - \frac{1}{2} +
O(e^{-2K\mu(r')}),
$$
{\it where} $r = \sqrt{t/(1+t)}.$

\bs

In \cite[Theorem 1.3]{32}, it was proved that, for $t \in (0, \infty)$ and
$K \in (1, \infty)$,
$$
16^{1-1/K} < \frac{\etK(t)}{t^{1/K}(1+t)^{K-1/K}} < 16^{K-1}
c^{2(1-1/K)}, \eqno(4.17)
$$
where $c = (1/8) \exp(2 \ck(1/\sqrt{2})^2/\pi) = 1.115\dots$.  The
proof of (4.17) given in \cite[Proof of Theorem 1.3]{32}, however,
is complicated and long.   Our next result improves (4.17), and
its proof is also much simpler.

\bs

{\bf 4.18. Theorem.}  (1)   {\it For each $t \in (0, \infty)$, the
function}
$$
f(K) \equiv 16^{1/K} \etK(t)/[t^{1/K} (1+t)^{K-1/K}]
$$
{\it is strictly increasing and convex from $[1,\infty)$ onto $[16,
\infty)$.  In particular, for $t \in (0, \infty)$ and} $K > 1$,
$$
\etK (t) > 16^{1-1/K} t^{1/K}(1+t)^{K-1/K}, \eqno(4.19)
$$
{\it and, for $t \in (0,\infty)$ and} $K \in (1,2)$,
$$
16^{1-1/K} < \frac{\etK(t)}{t^{1/K}(1+t)^{K-1/K}} < 16^{1-1/K}
\left\{1 + \left[\left(1+\sqrt{\frac{t}{1+t}}\right)^2-1\right](K-1)\right\}.
\eqno(4.20)
$$

(2) {\it For $t \in (0, \infty)$, let $r = \sqrt{t/(1+t)}, ~ A = A(r) =
r'e^{\mu(r')}$ and} $B = B(r) = r \exp \left(( 2/ \pi)
\ck(r)\ck'(r) - \mu(r')\right)$.  {\it Then the function}
$$
g(K) \equiv A^{-K} B^{1/K} \left[\frac{\etK(t)}{t^{1/K}
(1+t)^{K-1/K}} \right]^{1/2}
$$
{\it is strictly decreasing from $[1, \infty)$ onto $(1/4, B/A]$.  In
particular, for $t \in (0, \infty)$ and} $K \in (1, \infty)$,
$$
\begin{cases}
\ds{\frac{1}{16} A^{2K} B^{-2/K} < \frac{\etK(t)}{t^{1/K} (1+t)^{K-1/K}}
< A^{2(K-1)} B^{2(1-1/K)}} \\[.4cm]
\ds{ \qquad \qquad \qquad \le 16^{K-1} c^{2(1-1/K)}},
\end{cases} \eqno(4.21)
$$
{\it where $c$ is as in} (4.17).

\bs

{\bf Proof.}  Let $r = \sqrt{t/(1+t)}, ~ s = \varK(r)$.

For part (1), we rewrite $f$ as
$$
f(K) = \left( 4^{1/K} \frac{s}{s'} ~ \frac{r^{{\prime}\thinspace
K}}{r^{1/K}}\right)^2.
$$
By logarithmic differentiation and (1.21), we obtain
$$
\begin{cases}
\frac{1}{2} f'(K)  = f(K) \left[ \frac{1}{\mur} \ck'(s)^2 +
\frac{1}{K^2} \log \frac{r}{4} + \log r' \right]\\
\hspace{.56in} \equiv f(K) F_1(K).\end{cases} \eqno(4.22)
$$
Clearly, $F_1$ has the following limiting values:
$$
F_1(1) = m(r) + m(r') + \log \frac{rr'}{4}, ~ \lim_{K \to \infty}
F_1(K) = \mu(r') + \log r'. \eqno(4.23)
$$

Differentiation gives
$$
\frac{1}{2}K^3 F'_1(K) = F_2(K) \equiv \log \frac{4}{r} -
\frac{4}{\pi^2} \frac{\ck'(r)}{\ck(r)} \ck(s)^3 [\ce'(s) - s^2
\ck'(s)]. \eqno(4.24)
$$
It follows from \cite[Theorem 3.21 (1), (7)]{7}  that $F_2$ is strictly
increasing on $[1, \infty)$, and
\begin{align*}
\hfill \ds{F_2(1)} &= \ds{\left[ 1 - \frac{8}{\pi^3} \rp \ck(r)^3
\cdot
\frac{\mur}{\log(4/r)} \cdot \frac{\ce'(r) - r^2 \ck'(r)}{\rp} \right]
\mbox{log} \frac{4}{r}, }\\
\hfill \ds{\lim_{K \to \infty} F_2(K)} &= \ds{~\log \frac{4}{r}}.
\end{align*}

\medskip

 From \cite[Theorems 3.21 (1), (7) and 5.16 (2)]{7}, we
see that the function $F_2(1)/\log (4/r)$ is strictly increasing in
$r$ from $(0,1)$
onto $(0,1)$.  Hence, by (4.24), $F_1$ is strictly increasing on $[1,
\infty)$.  Since $F_1(1) > 0$ for all $r \in (0,1)$ \cite[Lemma
2.27]{32}, $F_1(K) > 0$ for all $K \ge 1$ and $r \in (0,1)$.  Hence,
$f(K)F_1(K)$ is positive and increasing on $[1, \infty)$ by (4.22),
and the monotonicity and convexity of $f$ follow from (4.22).

Clearly, $f(1) = 16$.  It follows from \cite[(2.11), p.  62]{20}  that
$$
\hspace{-.3in}  \lim_{K \to \infty} \frac{r^{\prime \thinspace K}}{s'}
= \lim_{K \to \infty} \frac{\exp[\mu(s')(1 +
(\log r')/\mu(r'))]}{\exp(\mu(s') + \log s')}
$$
$$
\qquad\qquad = \frac{1}{4} \lim_{s \to 1} \exp\left[
\frac{\mu(s')}{\mu(r')} (\mu(r') + \log r')\right] = \infty,
$$
and hence $\lim_{K \to \infty} f(K) = \infty$.

Inequality (4.19) is clear, while the upper bound in (4.20) follows
from the convexity of $f$ and (1.22).

For part (2), we write $g$ as
$$
g(K) = \frac{s}{s'} \left( \frac{r'}{A}\right)^K \left(\frac{B}{r}
\right)^{1/K}.
$$
Logarithmic differentiation gives
$$
\frac{K^2}{g(K)} g'(K) = g_1(K) \equiv \frac{4}{\pi^2} \mur \ck(s)^2
+ K^2 \log \frac{r'}{A} + \log \frac{r}{B} \eqno(4.25)
$$
and
$$
\frac{1}{2K} g'_1(K) = g_2(K) \equiv \frac{4}{\pi^2}
\frac{\ck(r)}{\ck'(r)} \ck'(s)^3 [\ce(s) - s^{{\prime}\thinspace 2}
\ck(s)] + \log \frac{r'}{A}.
$$
 From \cite[Theorem 3.21 (1), (7)]{7}  we see that $g_2$ is strictly
increasing on $[1, \infty)$ so that, for all $K \ge 1,$
$$
g_2(K) < \lim_{K \to \infty} g_2(K) = \mu(r') + \log(r'/A) = 0.
$$
Hence $g_1$ is strictly decreasing on $[1, \infty)$, so that
$$
g_1(K) < g_1(1) = \frac{2}{\pi} \ck(r) \ck'(r) + \log(rr'/(AB)) = 0
$$
for $K >1$.  Consequently, the monotonicity of $g$ follows from (4.25).

Clearly, $g(1) = B/A$.  From \cite[(2.11), p. 62]{20}, we get
$$
\lim_{K \to \infty} g(K) = \lim_{K \to \infty} \frac{s}{s'} \left(
\frac{r'}{A} \right)^K = \lim_{s \to 1} \frac{1}{s'e^{\mu(s')}} =
\frac{1}{4}.
$$
The first and second inequalities in (4.21) are clear.  The third
inequality in (4.21) holds if and only if $(A/4)^K (B/c) \le 1$ for all
$K > 1$ and $t > 0$.  This is true if and only if $AB/(4c) \le 1$ since
$A < 4$.  Now \cite[Lemma 2.27]{32}  implies that $AB \le 4c.  \qquad
\square$

\bs

{\bf 4.26. Remark.}  It follows from \cite[Theorem 5.13 (2)]{7}  and Lemma
2.2 (1) that $A^K B^{-1/K}$ is strictly increasing in $r$ from $(0,1)$
onto $(4^{-1/K}, 4^K)$.  Hence the lower bound in (4.21) tends to
$16^{K-1}$ as $r$ tends to 1, that is, as $t$ tends to $\infty$.
Consequently, for large $t$, the lower bound in (4.21) is better
than the lower bound $16^{1-1/K}$ given in (4.17).

\bs

{\bf 4.27. Theorem.}  {\it For $t \in (0, \infty)$, let} $r =
\sqrt{t/(1+t)}$.

(1)  {\it The function}
$$
f(K) \equiv \left\{ \etK(t) \exp\left(2\left(\frac{2}{\pi}
\ck(r) \ck'(r) + \log \frac{r'}{r}\right)\right)\right\}^{1/K}
$$
{\it is strictly
decreasing from $[1, \infty)$ onto $(e^{2 \mu(r')}, e^{4\ck(r)
\ck'(r)/\pi}]$.  In particular, for $t \in (0, \infty)$ and $K \in (1,
\infty),$}

$$
t \exp \left\{ 2\left[K \mu(r') - \frac{2}{\pi} \ck (r)\ck'(r)\right]\right\}
< \etK (t)
< t \exp\left\{ \frac{4}{\pi} (K-1) \ck (r)\ck'(r)\right\}.\eqno{(4.28)}
$$

(2)  {\it The function}
$$
g(K)\! \equiv\! \frac{1}{K\!-\!1}\! \log\! \left\{
\frac{\etK(t)}{t^{1/K} (1\!+\!t)^{K\!-\!1/K}}
 \exp \!\left[2\left(\frac{1}{K}\!-\!1\right)\!
\left(\frac{2}{\pi} \ck(r)
 \ck'(r)\! +\!\log(rr')\right) \right] \right\}
$$
{\it is strictly increasing from $(1, \infty)$ onto $(0,2[\mu(r') + \log
r']).$  In particular, for $t \in (0,\infty)$ and $K \in (1, \infty),$
$$
\begin{cases}
16^{1-1/K} ~~< \exp \left\{ 2\left( 1 - \frac{1}{K} \right) \left(
\frac{2}{\pi} \ck (r) \ck'(r) + \log(rr') \right) \right\} \\[.2cm]
\ds{ \qquad \qquad  < \frac{\etK(t)}{t^{1/K}(1+t)^{K-1/K}}} \\[.4cm]
\qquad \qquad < \exp \left\{ 2(K-1)  \right) \left[\mu (r') + \log r'
\right]   \\[.2cm]
\qquad \qquad   \left.  + 2 \left( 1 - \frac{1}{K} \right)
\left[\frac{2}{\pi} \ck (r) \ck'(r) + \log(rr') \right]
\right\}
\end{cases} \eqno(4.29)
$$
and
$$
e^{(K-1)[\log 2 + (a - \log 2)/K]} < \la(K) < e^{(K-1)[\pi + (a- \log
2)/K]},  \eqno(4.30)
$$
where} $a = (4/\pi) \ck(1/\sqrt{2})^2$.

\bs

{\bf Proof.}  Let $s = \varK(r)$.  Then $\etK(t) = (s/s')^2$.

(1)  Logarithmic differentiation gives
$$
\frac{K^2}{2f(K)} f'(K) = f_1(s) - f_1(r),
$$
where $f_1(x) = (2/ \pi) \ck(x) \ck'(x) + \log(x'/x)$.  Hence the
monotonicity of $f$ follows from Lemma 2.2 (2).

Clearly, $f(1) = e^{4 \ck(r) \ck'(r) /\pi}.$  By l'H\^opital's Rule,
(4.2), and \cite[(5.2)]{7}, we have
$$
\lim_{K \to \infty}  \etK(t)^{1/K} = \exp \left( \lim_{K \to
\infty} \frac{\log \etK(t)}{K} \right) = e^{2\mu(r')}.
$$
Hence $\lim_{K\to \infty} f(K) = e^{2 \mu(r')}$.

\medskip

(2) Let
$$
g_1(K) = \log(s/s') \!+\! K \log r'\! +\! ((1/K)\! - \!1)
 [ (2/\pi) \ck(r) \ck'(r)\! +\! \log r']\! - \!\log r
$$
and $g_2(K) = K-1$.  Then $g(K) = 2 g_1 (K)/g_2(K)$ and, by
differentiation,
$$
\frac{g'_1(K)}{g'_2(K)}
 = g_3 (K) \equiv \frac{1}{\mur} \ck'(s)^2 - \frac{1}{K^2} \left[
\frac{2}{\pi} \ck (r) \ck'(r) + \log r'\right]+ \log r' \eqno(4.31)
$$
and
$$
\begin{cases}
\ds{\frac{K^3}{2} g'_3 (K) = g_4(K)  \equiv \frac{2}{\pi} \ck(r) \ck'(r) +
\log r'} \\[.5cm]
\ds{\qquad \quad - \left( \frac{2}{\pi}\right)^2 \frac{\ck'(r)}{\ck(r)}
\ck(s)^3
\left[\ce'(s) - s^2 \ck'(s) \right].}  \end{cases} \eqno(4.32)
$$
 From \cite[Theorem 3.21 (1), (7)]{7}, we see that $g_4$ is strictly
increasing on $(1, \infty)$.  By (1.23) and (1.19),
$$
\begin{cases}
g_4(1)&=m(r')+\log r'+m(r)-\frac{4}{\pi^2} \ck(r)^2
\ck'(r)[\ce'(r)-r^2 \ck'(r)] \\
 &= m(r')+\log r' + \frac{4}{\pi^2} \ck(r)
\ck'(r) \{ \frac{\pi}{2} \rp -\ck(r)[\ce'(r)- r^2 \ck'(r)]\}\\
&= m(r') + \log r' + \frac{4}{\pi^2} \ck(r)
\ck'(r) \{ r^2 \ck(r) [\ck'(r)- \ce'(r)]\\
&\hspace{.7in}- \rp \ck'(r)[\ck(r) - \ce(r)]\} \\
&= m(r') +\log r'+ m(r)m(r') \left\{\ds{ \frac{\ck'(r)
-\ce'(r)}{\rp \ck'(r)}-\frac{\ck(r) - \ce(r)}{r^2 \ck(r)}}
\right\}.
\end{cases} \eqno(4.33)
$$
It follows from \cite[Corollary 3.9 (1)]{28}  and (4.33) that
$$
g_4(1) > 0 \quad \mbox{for} \quad r \in (0, 1/\sqrt{2}],
\eqno(4.34)
$$
and that
$$
g_4(1) > g_5(r) \equiv m(r') + \log r'- \frac{1}{2} m(r) m(r')
\eqno(4.35)
$$
for $r \in [1/\sqrt{2}, ~ 1).$

By \cite[Theorem 3.30 (1)]{7}  and Lemma 2.1 (1), we see that $g_5$ is
strictly increasing on $[1/\sqrt{2}, ~ 1)$.  Since
$$
g_5(1/\sqrt{2}) = \frac{1}{\pi} \ck(1/\sqrt{2})^2 - \frac{1}{2
\pi^2} \ck(1/\sqrt{2})^4 - \frac{1}{2} \log 2 = 0.148 \dots > 0,
$$
it follows from (4.35) that
$$
g_4(1) > 0 \quad \mbox{for} \quad r \in [1/\sqrt{2}, ~ 1).
\eqno(4.36)
$$

It follows from (4.34), (4.36), (4.32), and the monotonicity of $g_4$
that $g_3$ is strictly increasing on $(1, \infty)$.  Hence the
monotonicity of $g$ follows from (4.31) and \cite[Lemma 1.25]{7}.

The limiting values of $g$ follow from (4.31) and l'H\^opital's Rule.
The second and third inequalities in (4.29) are clear, while the first
follows from \cite[Lemma 2.27]{32}.  Taking $t=1$ in (4.29), we get (4.30).
$\qquad \square$

\bs

\section*{Acknowledgments.}

This work was completed during the second author's visits to Michigan
State University at East Lansing, MI, U.S.A., and to the University of
Helsinki, Finland.  His visits were supported by grants from the
Department of Mathematics at Michigan State University, the Finnish
Academy of Science and Letters, the Academy of Finland, and the Commission
on Development and Exchanges in Paris.  He wishes to express his thanks to
these institutions.

\bs
\bs

{\small

\vspace{.4in}

\begin{flushleft}
{\bf ANDERSON:} \\
Department of Mathematics \\
Michigan State University \\
East Lansing, MI 48824 \\
U. S. A.\\
e-mail: ~~{\tt anderson@math.msu.edu}\\[.6cm]

{\bf QIU:} \\
President \\
Zhejiang Sci-Tech University\\
Xiasha Higher Education Zone\\
Hangzhou 310038 \\
P. R. CHINA \\ [.6cm]

\noindent VUORINEN:\\
Department of Mathematics \\
FIN-20014 University of Turku \\
FINLAND\\
e-mail: ~~{\tt vuorinen@utu.fi} \\

\end{flushleft}

\end{document}